\documentclass[12pt]{article}%
\usepackage{graphics}
\usepackage{float}
\usepackage{amsmath}
\usepackage{amsfonts}
\usepackage{amssymb}
\usepackage{graphicx}%
\begin{document}

\title{Choreographies in the Discrete Nonlinear Schr\"odinger Equations}

\author{Renato Calleja\thanks{Instituto de Investigaciones en Matem\'{a}ticas
Aplicadas y en Sistemas, Universidad Nacional Aut\'{o}noma de M\'{e}xico,
calleja@mym.iimas.unam.mx}, Eusebius Doedel\thanks{Department of Computer
Science, Concordia University, Montreal, Canada, doedel@cs.concordia.ca},
Carlos Garc\'{\i}a-Azpeitia\thanks{Facultad de Ciencias, Universidad Nacional
Aut\'{o}noma de M\'{e}xico, cgazpe@ciencias.unam.mx}, 
Carlos L. Pando L.\thanks{Instituto de Física, Benemérita Universidad Autónoma de Puebla, Apartado Postal J-48, Puebla 72570, México}}

\maketitle

\begin{abstract}
We study periodic solutions of the discrete nonlinear Schr\"{o}dinger 
equation (DNLSE) that bifurcate from a symmetric polygonal relative 
equilibrium containing  $n$ sites. 
With specialized numerical continuation techniques and a varying physically
relevant parameter we can locate interesting orbits, including infinitely 
many choreographies. 
Many of the orbits that correspond to choreographies are stable, as 
indicated by Floquet multipliers that are extracted as part of the 
numerical continuation scheme, and as verified \textit{a posteriori} 
by simple numerical integration. 
We discuss the physical relevance and the implications of our results.
\end{abstract}

\section{Introduction}

In the last two decades there has been a growing interest in the study of
choreographic solutions ("choreographies") of the $n$-body and $n$-vortex 
problems.
Choreographies are periodic solutions where the bodies or the vortices 
follows the same path. 
The first choreography was discovered numerically for the case of three
bodies in \cite{Mo93}, and its existence was rigorously proved in 
\cite{ChMo00}. 
The term "choreography" was adopted for the $n$-body problem after the work
of Sim\'{o} \cite{Si00}.
Since then, variational methods \cite{FeTe04,BT04}, numerical minimization
\cite{ChGe02}, numerical continuation \cite{Ca}, and computer-assisted 
proofs \cite{KaSi07}, have been used to determine choreographies of the 
$n$-body problem; see also the references in these papers. 
For vortices in the plane, choreographies have been constructed for $3$ 
and $4$ vortices in \cite{Bo}. 
For $n$ vortices in a general bounded domain, choreographic solutions have 
been found close to a stagnation point of a vortex \cite{Ba0}, and close 
to the boundary of the domain \cite{Ba}.

In \cite{Ca}, \cite{GaIz13}, and \cite{ChFe08}, chorographies are found in
dense sets of Lyapunov families that arise from the stationary $n$-polygon 
of bodies in a rotating frame. 
The existence of these choreographies depends only on the symmetries of 
the equations in rotating coordinates, {\it i.e.}, these results can be 
extended to find choreographies of the $n$-vortex problem in the plane, 
in a disk, or in a surface of revolution. 
Such results can also be extended to the discrete nonlinear Schr\"{o}dinger
equation (DNLSE), which appears in the study of optical waveguide arrays,
and in Bose-Einstein condensates trapped in optical lattices, among many 
other applications \cite{Kr}.

While much research has been done on choreographies in the $n$-vortex and 
$n$-body problems, we are not aware of its extension to the study of the 
existence of choreographies in periodic lattices of $n$ sites, as modeled 
by the DNLSE.
The interest in $n$-body and $n$-vortex choreographies can perhaps be 
explained by the fact that variational methods are better suited for 
singular potentials. 
On the other hand, the continuation methods used in \cite{Ca} are very well
suited for locating choreographies for other symmetric potentials, such as 
in the DNLSE.

Evidence suggests that linearly stable choreographies in the $n$-body 
problem only exist for $n=3$, and in the $n$-vortex problem for 
$n=3,\cdots,7$, as a consequence of the fact that the polygonal relative 
equilibrium of $n$ bodies is unstable for $n\geq3$ \cite{Ca}, and for 
vortices for $n>7$ \cite{GaIz12}. 
On the other hand, the DNLSE has dense sets of stable choreographies, 
as a consequence of the fact that the DNLSE has stable polygonal relative 
equilibria for all $n$ \cite{Ga15}.
We note that boundary value continuation methods can determine unstable 
choreographies as easily as stable ones; a property not shared by most 
other techniques.
The aim of our paper is to investigate the existence of choreographies 
in the DNLSE using a boundary value continuation method. 
As illustrative examples we present a selection of choreographies,
most of them stable, in a periodic lattice for the case of $n=9$,
$17$, and $31$ sites.

The lack of previous work on detecting choreographies in the DNLSE may 
be related to the difficulties encountered in the measurement of phases 
in physical problems modeled by the DNLSE. 
However, such difficulties appear to be surmountable in nonlinear optics.
Within the context of nonlinear optics, several predictions of the DNLSE 
have been found experimentally in the last two decades \cite{Opt1},
such as the formation of discrete solitons in waveguide arrays. 
The use of suitable optical techniques, known as laser heterodyne 
measurements, to detect the field intensity and the phase \cite{Opt2,Opt3}, 
may open the door to experimental observation of stable choreographic 
solutions. 

In Section~2 we consider Lyapunov families of periodic orbits, and their 
relation to the existence of choreographies in a periodic lattice of 
Schr\"{o}dinger sites. 
In Section~3 we present methods to continue the Lyapunov families, and we 
exhibit a small selection of the many linearly stable choreographies that
we have determined.
In Section~4 we discuss our final remarks on choreographies.  
\section{Lyapunov families and choreographies}

In a rotating frame with frequency $\omega$, $q_{j}(t)=e^{i\omega t}%
u_{j}(t)$, the equation that describes the dynamics of a lattice of $n$ 
sites is given by the Hamiltonian system 
\begin{equation}
\dot{u}=\mathcal{J}\nabla H_{\omega}(u), \quad \text{where} \quad
H_{\omega}=\frac{1}{2}\sum_{j=1}^{n}\left(  \frac{1}{2}\left\vert
u_{j}\right\vert ^{4}+\omega\left\vert u_{j}\right\vert ^{2}-\left\vert
u_{j+1}-u_{j}\right\vert ^{2}\right)  . \label{Equations}
\end{equation}
The sites $u_{j}(t)\in\mathbb{C}$ satisfy periodic boundary conditions
$u_{j}(t)=u_{j+n}(t)$.
The equation of motion has explicit polygonal equilibrium solutions given 
by%
\begin{equation}
a_{j}=ae^{ij(\alpha\zeta)},\qquad\zeta=\frac{2\pi}{n}, \label{SW1}%
\end{equation}
for $\omega(a)=4\sin^{2}(m\zeta/2)-a^{2}$, $\alpha=1,\cdots,n$ and $a\in
\mathbb{R}^{+}$. 
These solutions correspond to relative equilibria in the non-rotating frame given by $q_{j}(t)$for $j=1,\cdots,n$.
The linearized Hamiltonian system at the polygonal equilibrium 
$\mathbf{a} =(a_{1},...,a_{n})$ is 
$$\dot{u}=\mathcal{J}D^{2}H_{\omega(a)}(\mathbf{a})u.$$
In \cite{Ga15} it is proved that the matrix 
$\mathcal{J}D^{2}H_{\omega (a)}(\mathbf{a})$ 
has a pair of imaginary eigenvalues $\pm i\nu_{k}$ for each
$k\in\{1,\cdots,n-1\}$ such that
\begin{equation}
\frac{a^{2}}{2\cos\alpha\zeta\sin^{2}k\pi/n}<1. \label{In}%
\end{equation}
It is also proved in \cite{Ga15} that for each $a$ 
and $k\in\{1,\cdots,n-1\}$ such that (\ref{In}) holds, 
the equilibrium $\mathbf{a}$ has a \emph{global family} 
of periodic solutions that arises from the normal modes 
of the polygonal equilibrium, and has the form
\begin{equation}
u_{j}(t)=e^{ij\alpha\zeta}u_{n}\left(  \nu t\pm jk\zeta\right)  \text{,}
\label{TW}%
\end{equation}
where $u_{n}(t)=a+\mathcal{O}(b)$ is a $2\pi$-periodic function,  
$\nu =\nu_{k}+\mathcal{O}(b)$ is the frequency and $b$ is a parameterization of the local family. 
The traveling waves (\ref{TW}) form two-dimensional families parametrized 
by the amplitude $a$ and the bifurcation parameter $b$. In the non-rotating frame, after rescaling time, these periodic solutions 
are traveling waves of the form%
\[
q_{j}(t)=e^{i\omega t/\nu}e^{ij\alpha\zeta}u_{n}
     \left(  t\pm jk\zeta\right) ,\qquad\omega=\omega(a)\text{.}%
\]
We say that a Lyapunov orbit is $\ell:m$ resonant if $\ell$ and $m$ are
relatively prime such that%
\[
\frac{\omega}{\nu}=\frac{\ell}{m}\text{,\qquad}k\ell-\alpha 
  m\in n\mathbb{Z}\text{.}%
\]
Such frequencies $\nu$ are dense in the set of real numbers.
For an $\ell:m$ resonant orbit, we have%
\[
q_{n}(t)=e^{it\omega/\nu}u_{n}(t)~.
\]
Since $e^{it\omega/\nu}=e^{it\ell/m}$ is $2\pi m$-periodic, the function
$q_{n}(t)=e^{it\omega/\nu}u_{n}(t)$ is $2\pi m$-periodic. 
Also, since%
\begin{equation}
 q_{n}(t-2\pi)=e^{i(t-2\pi)\omega/\nu}u_{n}(t-2\pi)=e^{-i2\pi\ell/m}%
 q_{n}(t),\label{symq}%
\end{equation}
the orbit of $q_{n}(t)$ is invariant under rotation by $2\pi/m$.
The solutions satisfy
\begin{align*}
q_{j}(t)  &  = e^{it(\omega/\nu)}u_{j}(t)
       =e^{it(\omega/\nu)}e^{ij\alpha\zeta }u_{n}(t+jk\zeta)\\
&  = e^{it(\omega/\nu)}e^{i\alpha j\zeta}e^{-i(\omega/\nu)
  \left( t+jk\zeta\right)  }q_{n}(t+jk\zeta)=e^{-ij
  \left(  (\omega/\nu)k-\alpha\right) \zeta}q_{n}(t+jk\zeta)\text{.}%
\end{align*}
Using the facts that $\omega/\nu=\ell/m$ and $\zeta=2\pi/n$, we have
\[
j\left(  k\frac{\omega}{\nu}-\alpha\right)  \zeta=2\pi j
   \left(  \frac{\ell k-\alpha m}{mn}\right)  =2\pi j \left( \frac{r}{m} \right) ,
\]
with $r=(k\ell-\alpha m)/n\in\mathbb{Z}$ by assumption. 
Since $\ell$ and $m$ are relatively prime we can find $\ell^{\ast}$, 
the $m$-modular inverse of $\ell$. 
Since $\ell\ell^{\ast}=1$ mod $m$, it follows from the symmetry
(\ref{symq}) that
\[
q_{n}(t-2\pi jr\ell^{\ast})=e^{-i2\pi j(r/m)}q_{n}(t).
\]
Then%
\begin{equation}
q_{j}(t)=e^{-i2\pi j(r/m)}q_{n}(t+jk\zeta)
   = q_{n}(t+j(k-rn\ell^{\ast})\zeta).
\end{equation}
Thus in the non-rotating frame, an $\ell:m$ resonant Lyapunov orbit 
is a choreography satisfying%
\[
q_{j}(t)=q_{n}(t+j\tilde{k}\zeta)\text{,}%
\]
where 
$\tilde{k}=k-(k\ell-\alpha m)\ell^{\ast}$ 
with $\ell^{\ast}$ the $m$-modular inverse of $\ell$. 
The period of the choreography is $m~T_{\ell :m}$ with%
\[
T_{\ell:m}=\frac{2\pi}{\nu}=2\pi\omega(a)\left(  \frac{\ell}{m}\right)  .
\]
The choreography is symmetric under rotation by $2\pi/m$, and it winds
around a center $\ell$ times.
\section{Computational Results}
\label{sec:2}
We have computed families of periodic solutions that arise directly or
indirectly from the circular polygonal relative equilibrium of the DNLSE.
These families are computed by numerical continuation using boundary value 
formulations. 
In this article we present numerical results for several choices of the 
number of sites $n$ in the DNLSE, namely for $n=9$, $n=17$, and $n=31$,
with various values of the amplitude parameter $a$.
In our boundary value setting the DNLS differential equations are 
formulated as
\begin{align*}
u_{k}^{\prime}(t)  &  =-iT\left(  u_{k-1}-2u_{k}+u_{k+1}+\left\vert
u_{k}\right\vert ^{2}u_{k}+\omega u_{k}\right) \\
&  +p_{1}(4u_{k}-4u_{k}^{3}-2\bar{u}_{k+1})+p_{2}(u_{k+1}-\bar{u}_{k-1}),
\end{align*}
where $u_{k}(t)=x_{k}(t)+iy_{k}(t)$ for $k=1,2,\cdots,n$ and
\[
u_{0}(t)\equiv u_{n}(t)\quad,\quad u_{n+1}(t)\equiv u_{1}(t).
\]
Here $T=2\pi/\nu$ is the period of a periodic orbit, so that the scaled time
variable $t$ takes values in the fixed time interval $[0,1]$. The parameters
$p_{1}$ and $p_{2}$ are \textit{unfolding parameters} that are necessary to
take care of invariances related to the presence of two conserved quantities.
The parameters $p_{1}$ and $p_{2}$ are always part of the unknowns solved in
the Newton iterations. However, upon converge their values are zero up to
numerical accuracy.
The boundary conditions that we impose always include the periodicity
conditions
\[
u_{k}(1)-u_{k}(0)=0,\quad k=1,2,\cdots,n.
\]
Additional constraints can be used to fix certain quantities along 
solution families, provided other appropriate parameters are allowed 
to vary. 
In particular, we can fix the $y$-coordinate of the $n$th site at time 
$t=0$, {\it i.e.},%
\[
y_{n}(0)=0.
\]
This constraint can be viewed as a phase condition that is sometimes more
convenient than an integral phase condition of the type mentioned below. 
It can also be useful to fix the $x$-coordinate of the $n$th site 
at time $t=0$, which is accomplished by adding the boundary condition
\[
x_{n}(0)-x_{n}^{0}=0,
\]
where $x_{n}^{0}$ is a parameter that can be kept fixed.
For convenience, constraints of this form can also be used to keep track 
of such quantities. 
For example, the parameter $x_{n}^{0}$, when free to vary, trivially keeps 
track of $x_{n}(0)$. 
Such constraints can also fix (or to keep track of) one of the conserved 
quantities $E$ or $A$, or the resonance ratio $T/T_{0}$, 
where $T_{0}=2\pi/\omega$ is the period of the rotating frame. 
This is accomplished by adding one of the following constraints:
\[
\left\vert u_{k}-u_{k+1}\right\vert ^{2}-\left\vert u_{k}\right\vert
^{4}-\omega\left\vert u_{k}\right\vert ^{2}-E=0~,~~\mbox{where}\quad
\omega=4\sin^{2}(\pi/n)-a^{2},
\]%
\[
\sum_{k=1}^{n}\left\vert u_{k}\right\vert ^{2}-A=0,\quad\mbox{or}\quad
T/T_{0}-r=0.
\]
The boundary value formulation can also contain integral constraints, 
such as the phase condition
\[
\int_{0}^{1}x_{n}(t)~\tilde{x}_{n}^{\prime}(t)
          +y_{n}(t)~\tilde{y}_{n}^{\prime }(t)~dt=0,
\]
which here is applied to the $n$th site only, and where 
$(\tilde{x}_{n}^{\prime}(t),\tilde{y}_{n}^{\prime}(t))$ represents 
the time-derivative of a reference solution, which typically is the 
preceding solution in the numerical continuation process. 
Another integral constraint sets the average $y$-coordinate 
of the $n$th site to zero, namely,
\[
\int_{0}^{1}y_{n}(t)~dt=0.
\]
The purpose of this constraint is to remove the rotational invariance of
periodic solutions. 
There are more general integral constraints for fixing the phase and for 
removing invariances. 
However the ones listed above are simple, and appropriate in the current 
context.

\begin{figure}[ptb]
\begin{center}
\resizebox{16.0cm}{!}{  
\includegraphics{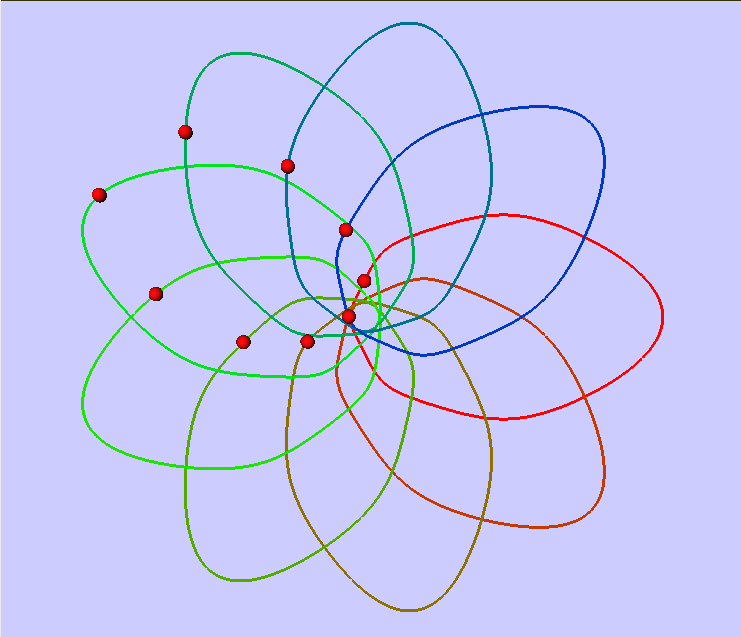}~
\includegraphics{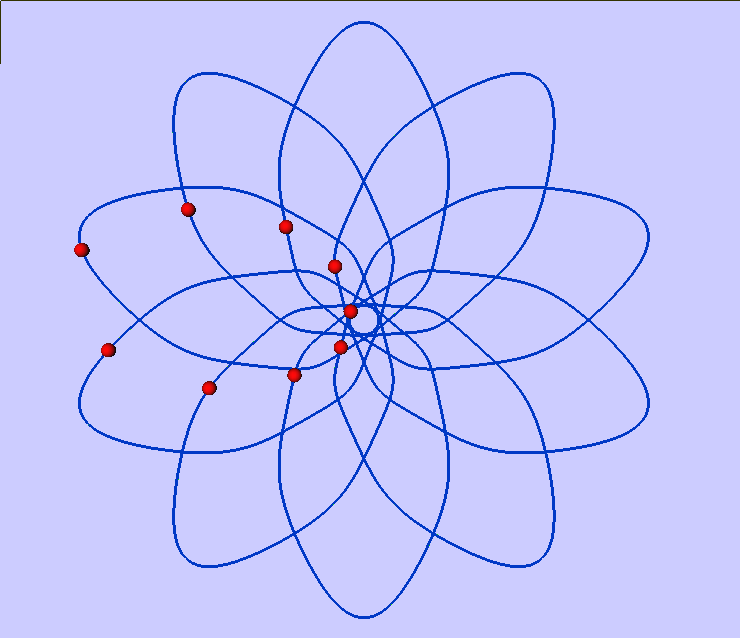} }
\end{center}
\par
\vskip-.90cm\noindent
\par
\begin{center}
\resizebox{16.0cm}{!}{  
\includegraphics{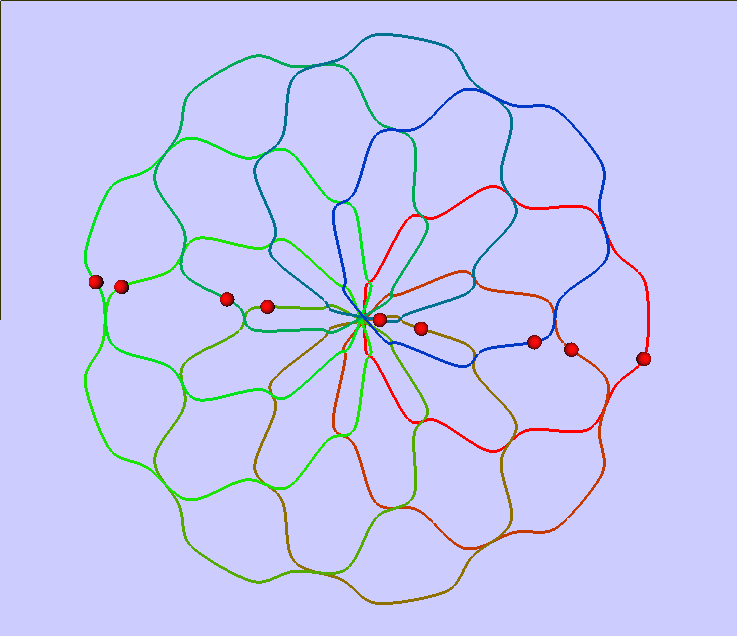}~
\includegraphics{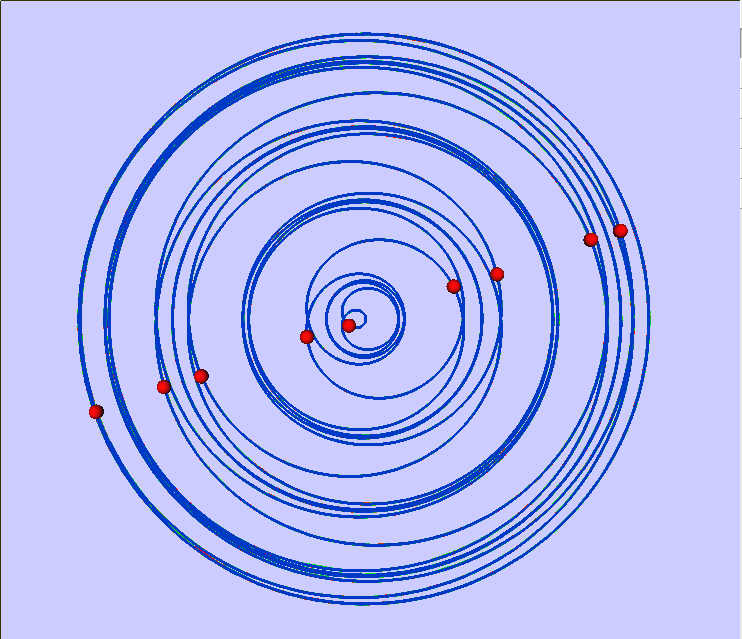} }
\end{center}
\par
\vskip-.90cm\noindent
\par
\begin{center}
\resizebox{16.0cm}{!}{  
\includegraphics{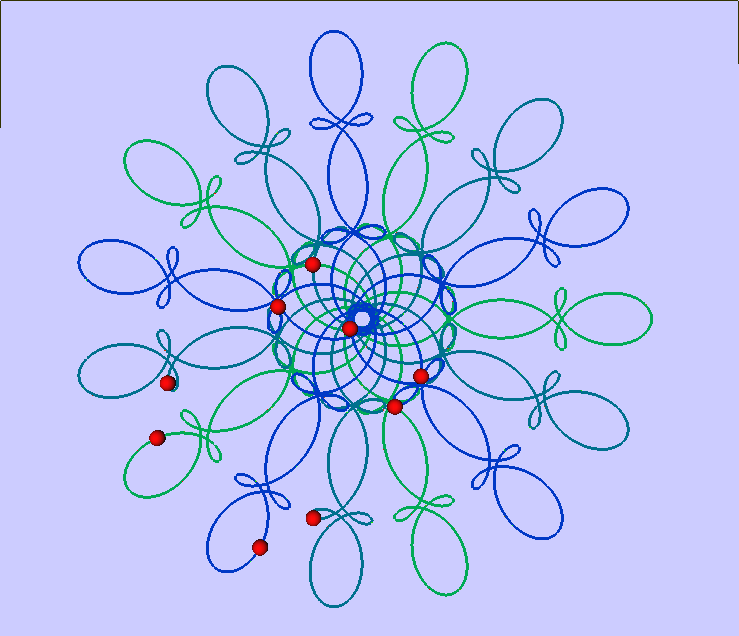}~
\includegraphics{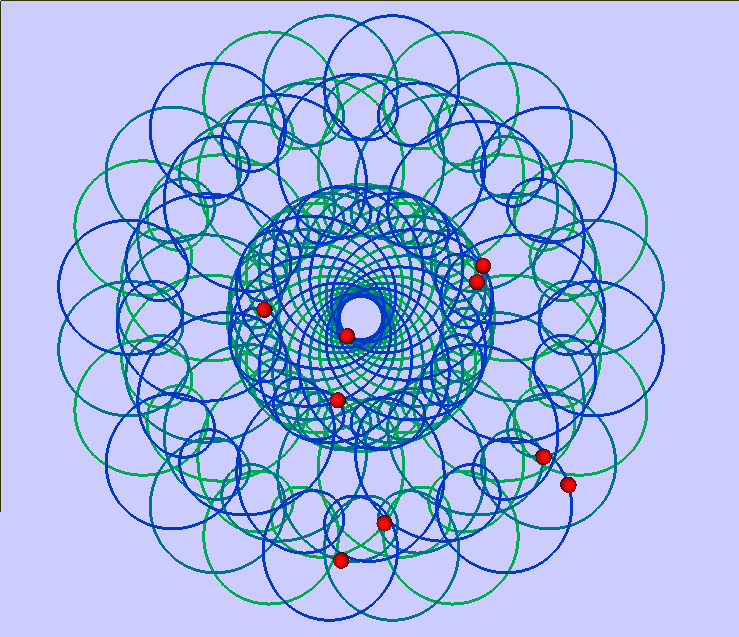} }
\end{center}
\par
\vskip-.5cm\noindent
\par
\caption{ 
All solutions in this figure are for $n=9$ and, to numerical 
accuracy, linearly stable.
Data are given in Table~1.
Top-Left:
A periodic solution in the rotating frame, of resonance 1:10. 
Top-Right:
The corresponding periodic solution in the non-rotating frame, 
where it is a choreography. 
Center-Left: 
A periodic solution in the rotating frame, of resonance 23:1. 
Center-Right:
The corresponding choreography.
Bottom-Left:
A partial choreography in the non-rotating frame, of resonance 2:5.
Bottom-Right: 
A partial choreography in the non-rotating frame of resonance 5:8.
}
\label{fig1}
\end{figure}

\begin{figure}[ptb]
\begin{center}
\resizebox{17.0cm}{!}{  
\includegraphics{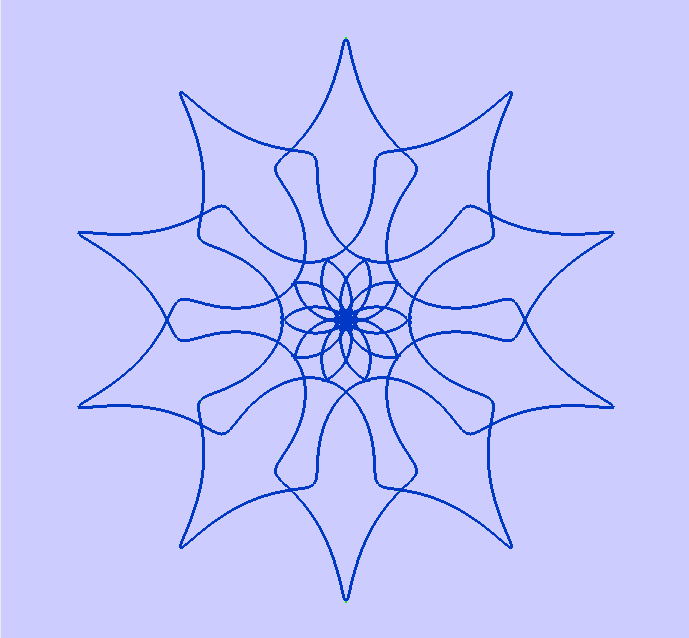}~~
\includegraphics{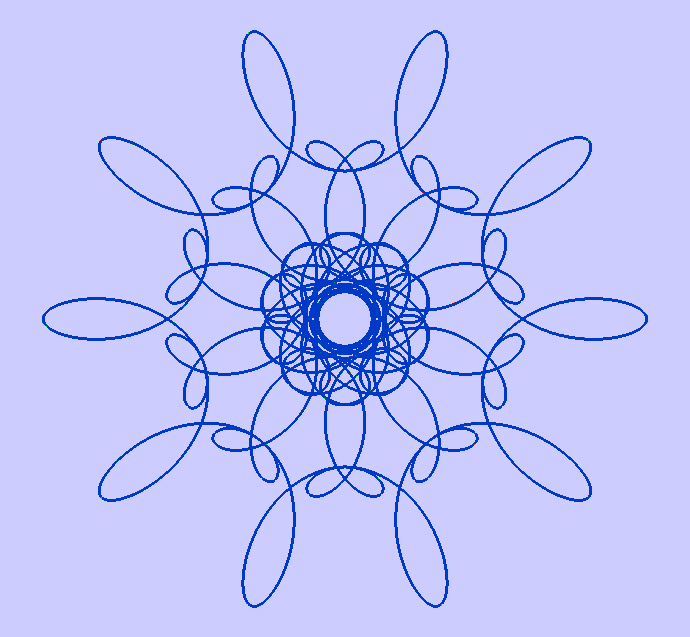}~~
\includegraphics{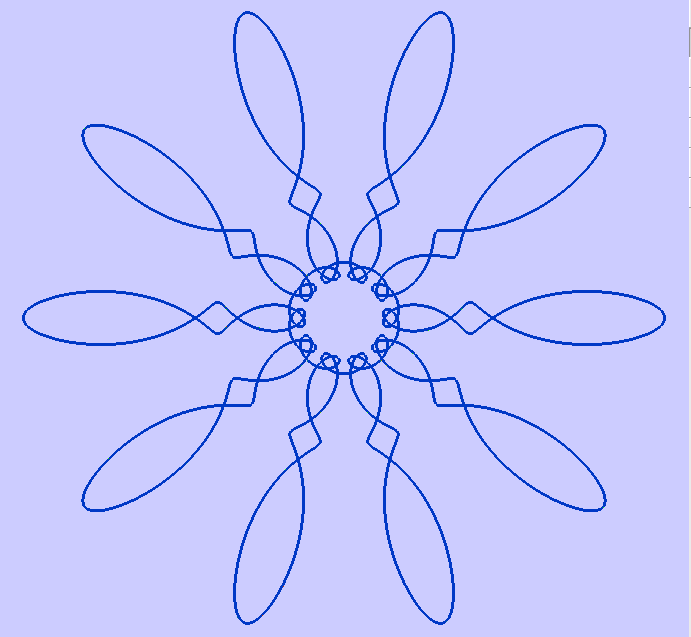} }
\end{center}
\par
\vskip-.90cm\noindent
\par
\begin{center}
\resizebox{17.0cm}{!}{  
\includegraphics{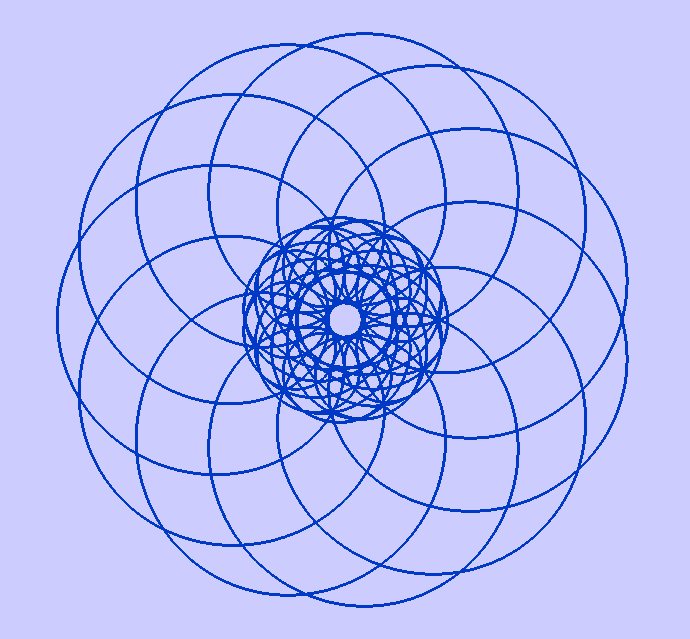}~
\includegraphics{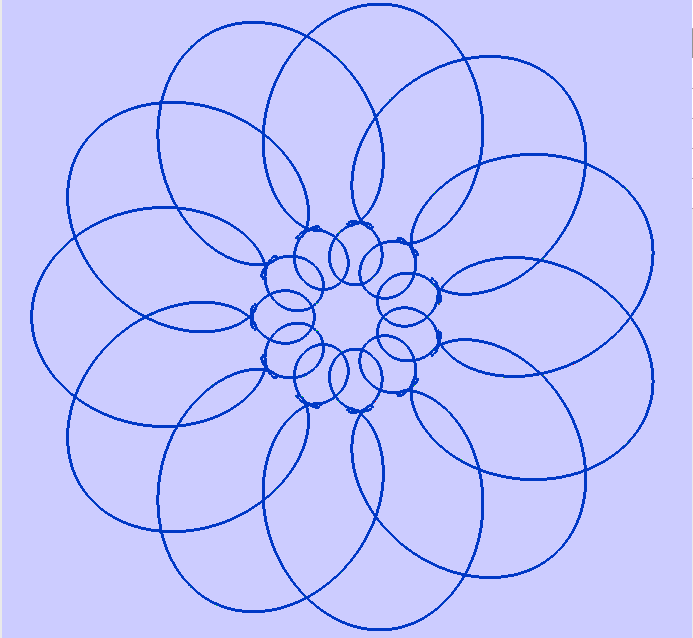}~
\includegraphics{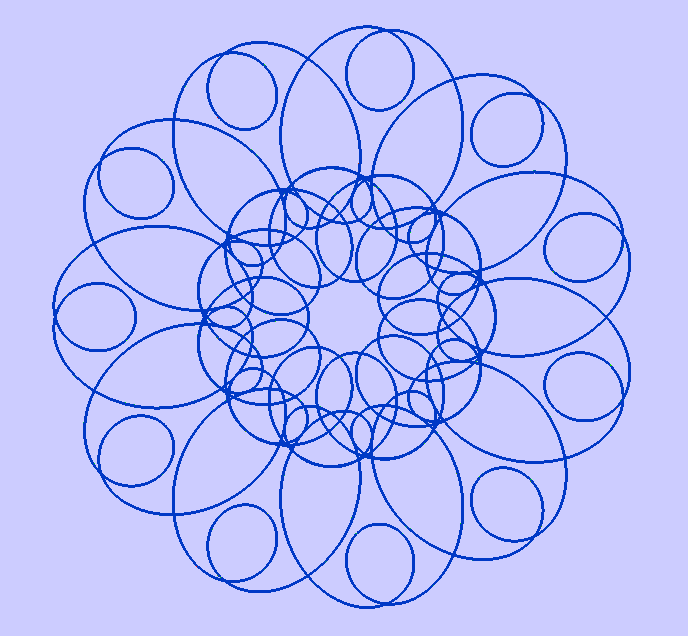} }
\end{center}
\par
\vskip-.90cm\noindent
\par
\begin{center}
\resizebox{17.0cm}{!}{  
\includegraphics{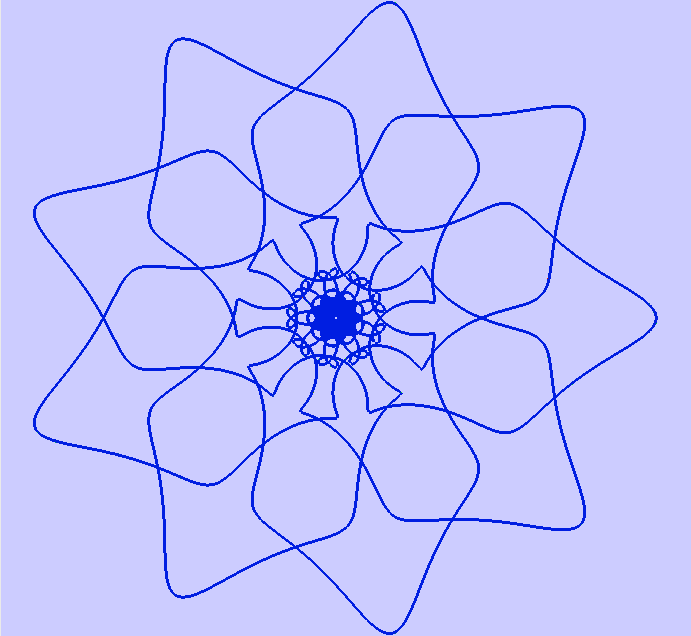}~
\includegraphics{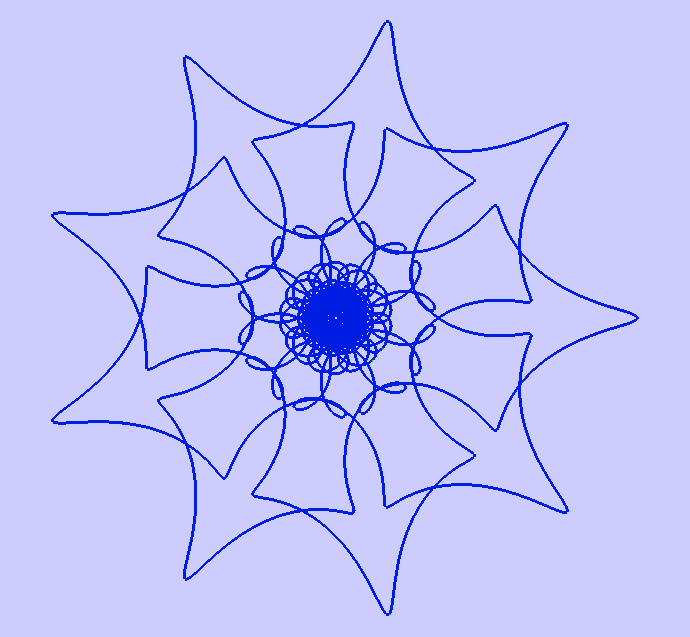}~
\includegraphics{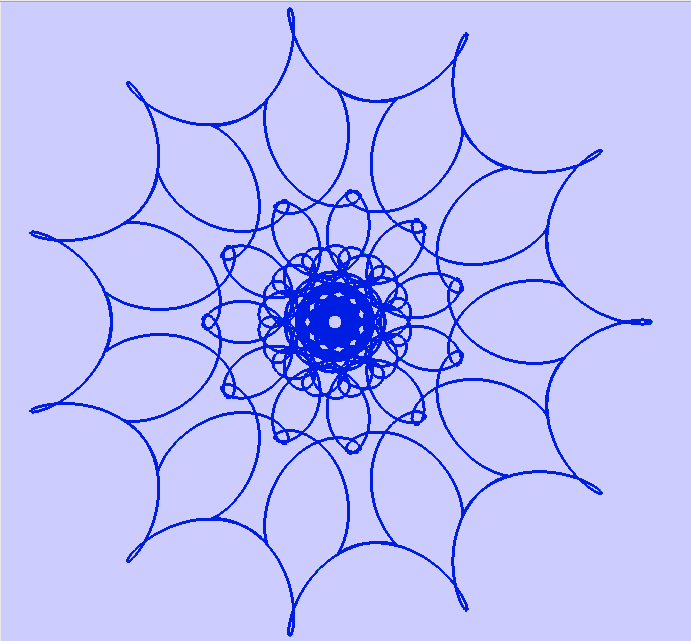} }
\end{center}
\par
\vskip-.90cm\noindent
\par
\begin{center}
\resizebox{17.0cm}{!}{  
\includegraphics{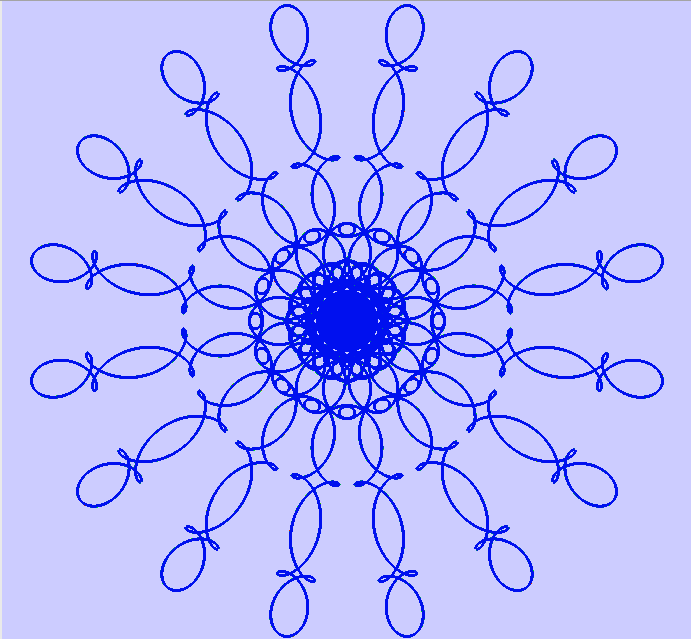}~
\includegraphics{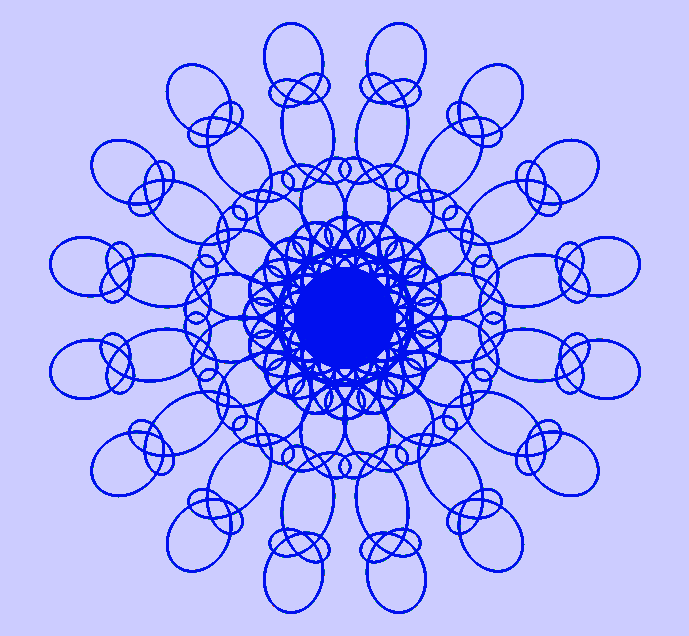}~
\includegraphics{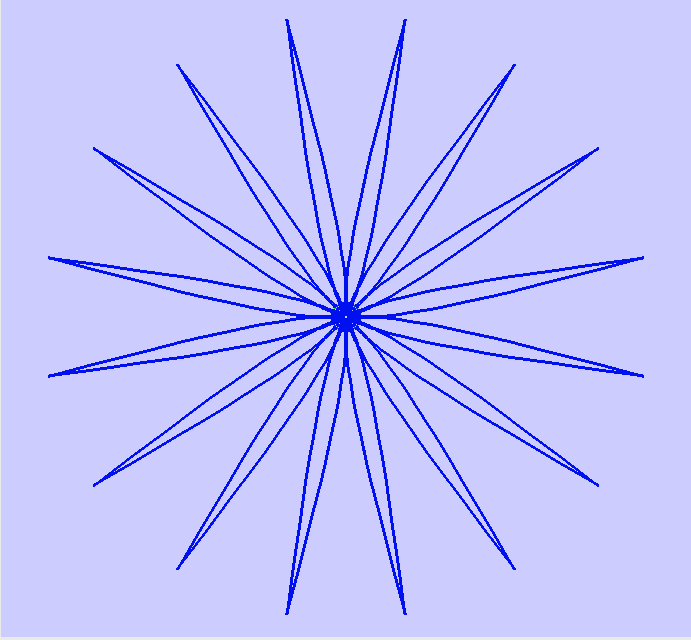} }
\end{center}
\par
\vskip-.5cm\noindent
\par
\caption{ 
The choreographies in the first two rows are for $n=9$,
the choreographies in the third row are for $n=17$, 
and the choreographies in the last row are for $n=31$.
All but two choreographies in this Figure are stable.
The two unstable choreographies are on the right in rows 1 and 2.
}
\label{fig2}%
\end{figure}

We now briefly outline the computational procedure that we have used to 
locate the periodic orbits shown in Figure~\ref{fig1} and in
Figure~\ref{fig2}, for which corresponding data is given in Table~1.
As a starting procedure we follow a family of periodic solutions that 
emanates from the polygonal equilibrium, namely a family that arises 
from a conjugate pair of purely imaginary eigenvalues of the equilibrium.
This starting procedure is relatively standard, and essentially the same 
as the one used for Hopf bifurcation. 
Initially only a small portion of the periodic solution family is computed;
in fact only a few continuation steps are taken. 
A minor adjustment of the standard starting procedure ensures that the 
small-amplitude starting solution satisfies in particular the conditions 
$y_{n}(0)=0$ and $\int_{0}^{1} y_{n}(t)~dt = 0$.
Keeping $y_{n}(0)=0$, the small-amplitude starting solution is followed 
until $x_{n}(0)$ reaches a specified target value, for which we have used 
values such as $x_{n}^{0}=-0.04$, $x_{n}^{0}=0.005$, and $x_{n}^{0}=0.0001$
(see Table~1).
The free continuation parameters in this step include $T$, $x_{n}^{0}$,
$p_{1}$, and $p_{2}$. 
The resulting periodic orbit has the property that all nine solution 
components pass near the origin in the complex plane. 
The subsequent main computational step then consists of keeping $x_{n}^{0}$ 
fixed at the value $x_{n}^{0}$, while allowing the amplitude parameter 
$a$ to vary. 
Specifically, the free continuation parameters now include 
$T$, $a$, $p_{1}$, and $p_{2}$. 
In the quest for locating interesting, stable  periodic solutions of the 
DNLSE, there are multiple variations on the continuation scheme outlined 
above. 
For example, one can start from a selected periodic solution from the main
computational step, now keeping the conserved quantity $H$ fixed. 
Yet another variation that we use is to follow periodic solutions found 
in the main computational step above, keeping the resonance ratio 
$T/T_0$ fixed. Here $T$ is the period of the periodic orbit and $T_{0}$ 
is the period of the rotating frame.
In fact, along all families of periodic solutions we monitor the value 
of the ratio $T/T_{0}$.
Specifically we are interested in rational values of $T/T_{0}$, for which 
the orbits correspond to a choreographies in the non-rotating frame.
As discussed in Section~2, there is a countably infinite number of such
choreographies along families of periodic orbits, provided that the 
the orbits possess certain symmetries, and provided the period $T$ is 
not constant.
Moreover, as already mentioned above, such choreographies can subsequently
be continued with varying amplitude parameter $a$, while keeping the ratio
$T/T_{0}$ fixed at a choreographic value.
Thus the number of choreographies then becomes in fact uncountably infinite.

Two choreographies for $n=9$ are shown in Figure~\ref{fig1}, namely 
in the top-right panel and center-right panel of that Figure.
In addition, all twelve orbits shown in Figure~\ref{fig2} correspond to
choreographies.
Specifically, in Figure~\ref{fig1}, the top-left panel shows a resonant 
orbit in the rotating frame.
The coloring of this orbit is according to its nine components, 
and evidently it consists of nine separate closed curves.
The top-right panel of Figure~\ref{fig1} shows the same orbit in the
non-rotating frame, where all components follow a single curve, {\it i.e.},
the orbit is a choreography.
The coloring of the choreography is also according to its nine components.
However, since all components follow the same curve, the color of this
curve changes gradually to a uniform final color, as one complete orbit
is traversed.

Similarly, the center-left panel of Figure~\ref{fig1} shows an orbit in
the rotating frame, while the center-right panel shows the same orbit in
the non-rotating frame, where it evidently corresponds to a choreography.
The bottom panels of Figure~\ref{fig1} show two resonant orbits in the
non-rotating frame.
As the coloring indicates, neither of these two orbits is a choreography.
However, both can be designated as a {\it partial choreography}, since
each consists of three separate closed curves, and each of these three
curves is traversed by three components.
Finally, Figure~\ref{fig2} shows a selection of twelve orbits in the
non-rotating frame, each of which is a choreography.
These choreographies are visually appealing, and particularly interesting 
to watch in animations.

\section{Conclusions and discussion}
\label{sec:3}
The DNLS equations with $n$ sites, and with periodic boundary conditions, 
have Lyapunov families of periodic solutions that arise from polygonal 
equilibria in the rotating frame. These Lyapunov families can be 
parameterized by the rotational frequency $\omega$ and the frequency 
$\nu$ of the Lyapunov orbit. 
When the ratio of the freqencies $\omega$ and $\nu$ is $l:m$ resonant, 
{\it i.e.}, when $\omega/\nu=l/m$, then the Lyapunov orbit corresponds 
to a choreography that is symmetric with respect rotations by $2\pi/m$, 
and that has winding number $l$. 
For fixed $\omega$ this resonance condition is satisfied for a dense set 
of rational frequencies $\nu$.
Thus if the frequency range of $\nu$ along a Lyapunov family contains 
an interval, then there is an infinite number of Lyapunov orbits that 
correspond to choreographies.

In this paper a robust and highly accurate boundary value technique with 
adaptive meshes has been used to continue the Lyapunov families. 
The presence of two conserved quantities, namely the amplitude $A$ and the 
energy $E$, is dealt with by using two unfolding parameters. 
The formulation also allows continuation of solution families with fixed 
$E$, $A$, or $\omega/\nu$. 
For example, by fixing $\omega/\nu=l/m$ one can compute a continuum of
choreographies, each of which is symmetric with respect to rotation by
$2\pi/m$, and has winding number $l$.  We have included a small sample 
of the infinitely many stable choreographies that can be computed 
in this manner.

In principle, physical observation of stable choreographies appears to be
possible with heterodyne optical techniques.
Such techniques have been used experimentally to record both the phase 
and the amplitude of a coherent optical signal source, such as that at 
the end of a waveguide array, as modeled by the DNLSE. 
The basic idea is to study the optical field composed of a reference and 
a source field.
This optical technique has allowed the confirmation of predictions from
the Lorenz model, which describes to a good degree the dynamics of the 
$NH_{3}$-laser \cite{Opt2}. 
Similarly, heterodyne optical techniques have been used to study optical 
fields in two-dimensional light sources \cite{Opt3}, which are precisely 
the geometries that support stable choreographies.
\clearpage
\clearpage
\clearpage
\vskip.0cm\noindent
\begin{table}[H]
\begin{center}
\begin{large}
\begin{tabular}{|c|c|c||c|c|c|c|c|c|c|}
\hline
      &   &   &   &       &          &         &          &         &   \cr
Figure&Row-&Orbit-&$n$&$T/T_0$& a    & $T$     &   $T_0$  &$x_n(0)$ &S/U\cr
      &Col.&Label&&        &          &         &          &         &   \cr
\hline
      &1-1& 1 & 9 & 1:10  & 0.651774 & 14.5773 &  145.773 & -0.04   & S \cr
      &1-2& 2 & 9 & 1:10  & 0.651774 & 14.5773 &  145.773 & -0.04   & S \cr
   1  &2-1& 3 & 9 & 23:1  & 0.657102 & 4000.00 &  173.913 &  0.005  & S \cr
      &2-2& 4 & 9 & 23:1  & 0.657102 & 4000.00 &  173.913 &  0.005  & S \cr
      &3-1& 5 & 9 & 2:5   & 0.520316 & 12.7459 &  31.8649 & -0.04   & S \cr
      &3-2& 6 & 9 & 5:8   & 0.396319 & 12.6334 &  20.2134 & -0.04   & S \cr
\hline
      &1-1& 7 & 9 & 1:10  & 0.647930 & 13.0635 &  130.635 & -0.04   & S \cr
      &1-2& 8 & 9 & 1:10  & 0.646671 & 12.6423 &  126.353 & -0.04   & S \cr
      &1-3& 9 & 9 & 1:10  & 0.627791 & 8.51610 &  85.1505 & -0.04   & U \cr
      &2-1&10 & 9 & 2:11  & 0.510285 & 5.50498 &  30.2774 & -0.04   & S \cr
      &2-2&11 & 9 & 2:11  & 0.531986 & 6.17839 &  33.9811 & -0.04   & S \cr
   2  &2-3&12 & 9 & 2:11  & 0.565906 & 7.73660 &  42.5513 & -0.04   & U \cr
      &3-1&13 &17 & -8:9  & 0.576588 & 28.2933 & -31.8299 &  0.005  & S \cr
      &3-2&14 &17 & -8:9  & 0.578005 & 28.0608 & -31.5684 &  0.005  & S \cr
      &3-3&15 &17 & -6:11 & 0.505528 & 28.4406 & -52.1411 &  0.005  & S \cr
      &4-1&16 &31 & -15:16& 0.421561 & 43.0673 & -45.9385 &  0.0001 & S \cr
      &4-2&17 &31 & -15:16& 0.421549 & 43.0706 & -45.9420 &  0.0001 & S \cr
      &4-3&18 &31 & -15:16& 0.420348 & 43.3913 & -46.2840 &  0.0001 & S \cr
\hline
\end{tabular}
\end{large}
\caption{
Data for the orbits in Figures~1 and 2.
Here $n$ is the number of sites,
$a$ the amplitude parameter, 
$T$ the period of the orbit,
$T_0$ the period of the rotating frame,
$x_n(0)$ the $x$-component of the $n$th site at time zero,
"S" stands for "Stable or almost stable", and
"U" stands for "Unstable".
}
\end{center}
\end{table}

\section*{Acknowledgments}
This work was supported by NSERC (Canada), BUAP and CONACYT (M\'{e}xico).

\end{document}